\documentclass{conm-p-l}
\usepackage{amsmath,amsfonts,amsthm,amssymb,eucal,amsxtra}

\newcommand{\Qp}{\mathbb Q_p}
\newcommand{\Cp}{\mathbb C_p}

\newcommand{\LA}{\mathcal L_A}
\newcommand{\LP}{\mathcal L_P}
\newcommand{\B}{\mathcal B}
\newcommand{\MA}{\mathcal M(\LA)}

\newcommand{\AAA}{\mathfrak A}

\newcommand{\AR}{\widetilde{\mathfrak A}}
\newcommand{\GD}{\widehat{G}}
\newcommand{\II}{\mathbf I}
\newcommand{\JJ}{\mathbf J}
\newcommand{\RR}{\mathbf R}
\copyrightinfo{2012}{American Mathematical Society}

\begin{document}
\newtheorem{lem}{Lemma}[section]
\newtheorem{teo}{Theorem}[section]
\newtheorem{cor}{Corollary}[section]
\numberwithin{equation}{section}
\pagestyle{plain}
\title{On Some Classes of Non-Archimedean Operator Algebras}
\author{Anatoly N. Kochubei}
\address{Institute of Mathematics,National Academy of Sciences of Ukraine,Tereshchen\-kivska 3, Kiev, 01601 Ukraine}
\thanks{This work was supported in part by Grant No. 01-01-12 of the National Academy of Sciences of Ukraine under the program of joint Ukrainian-Russian projects.}
\email{kochubei@i.com.ua}
\date{}
\subjclass[2010]{Primary: 47S10. Secondary: 47L10; 47L65; 16S99.}
\bigskip
\begin{abstract}
We study some classes of algebras of operators on non-Archime\-dean Banach spaces.  In particular, we propose a non-Archimedean version of the crossed product construction.
\end{abstract}
\maketitle

\section{Introduction}

In this paper we consider some classes of algebras of bounded linear operators on Banach spaces over non-Archimedean fields. Firstly, in terms of a Banach algebra generated by a single operator, we define its normality property. In Section 2, we prove the strong normality (the definitions are given below) of the orthoprojections (in the non-Archimedean sense) and multiplication operators on function spaces. We also consider commutative algebras of normal operators. For other approaches in this direction see \cite{LD}.

Secondly (Section 3) we propose a possible way based on the notion of a Baer ring \cite{Kap}) to develop a counterpart of the notion of a von Neumann algebra. The idea is based on a reduction procedure. Given a non-Archimedean operator algebra $\AAA$, this procedure produces a ring $\AR$ of operators on a vector space over the residue field (see \cite{BGR} for various applications of this idea). For algebras $\AAA$ whose reductions belong to the class of Baer rings, this leads to a kind of classification.
Short of any general theory of that kind, we consider in detail a non-Archimedean version of the crossed product construction, one of the main methods of constructing von Neumann algebras in the classical case (see the original paper \cite{MN} by Murray and von Neumann; for a modern exposition see \cite{Ta}). This results (Section 4) in a class of non-trivial non-Archimedean factors (algebras with a trivial center) or algebras close to factors corresponding, through the reduction procedure, to type I Baer rings.

\section{Normal Operators}

In this section we recall the notions from \cite{K1} and provide some corrections and additional material. We will not explain the basic notions of non-Archimedean analysis; see \cite{BGR,PS,Rob,Rooij,Sch}.

{\bf 2.1.} {\it Spectral theorem.} Let $A$ be a bounded linear operator on a Banach space $\mathcal B$ over a complete non-Archimedean valued field $K$ with a nontrivial valuation; $|\cdot |$ will denote the absolute value in $K$, $O$ is the ring of integers in $K$. We denote by $\| \cdot \|$ both the norm in $\mathcal B$ and the operator norm $\|A\|=\sup\limits_{x\ne 0}\dfrac{\|Ax\|}{\|x\|}$ (see \cite{NBB} for some subtleties regarding the operator norm in the non-Archimedean case). Below we assume that $K$ is algebraically closed. We will denote by $\widetilde{K}$ the residue field of $K$.

Denote by $\LA$ the commutative Banach algebra generated by $A$ and the unit operator $I$. $\LA$ is the closure of the algebra $K[A]$ of polynomials in $A$, with respect to the norm of operators; thus $\LA$ is a Banach subalgebra of the algebra $L(\B)$ of all bounded linear operators. Elements $\lambda \in K$ are identified with the operators $\lambda I$.

The spectrum $\MA$ of the algebra $\LA$ is defined (see \cite{Ber}) as the set of all bounded multiplicative seminorms $|\cdot |$ on $\LA$. $\MA$ is provided with the weakest topology, with respect to which all functions $|\cdot | \mapsto |B|$, $B\in \LA$, are continuous.
In this topology, it is a nonempty Hausdorff compact topological space. If the algebra $\LA$ is uniform, that is $\|T^2\|=\|T\|^2$ for any $T\in \LA$, and all the characters take their values in $K$, then \cite{Ber} the space $\MA$ is totally disconnected, and $\LA$ is isomorphic to the algebra $C(\MA ,K)$ of continuous functions on $\MA$ with values in $K$. In this case the above isomorphism transforms the characteristic functions $\eta_A$ of nonempty open-closed subsets $\Lambda \subset \MA$ into idempotent operators $E(\Lambda)\in \LA$, $\|E(\Lambda)\|=1$. These operators form a finitely additive norm-bounded projection-valued measure on the Boolean algebra of open-closed sets, with the non-Archimedean orthogonality property
$$
\|f\| =\sup\limits_\Lambda \|E(\Lambda )f\|,\quad f \in \B.
$$

An operator with the above properties is called {\it normal}. It is called {\it strongly normal}, if
its spectrum $\sigma (A)$ is a nonempty totally disconnected compact
subset of $K$, and $\MA =\sigma (A)$. For a strongly normal operator $A$, we have
the spectral decomposition
$$
A=\int\limits_{\sigma (A)}\lambda \,E(d\lambda ).
$$
For any $\varphi\in C(\sigma (A),K)$ we can define the operator
$$
\varphi (A)=\int\limits_{\sigma (A)}\varphi (\lambda )\,E(d\lambda ).
$$
The operator $\varphi (A)$ is strongly normal \cite{K3}.

In order to check the normality of an operator $A$, one has to establish the equality
\begin{equation}
\|[q(A)]^2\|=\|q(A)\|^2
\end{equation}
for any polynomial $q$ over $K$. Unfortunately, a simple sufficient condition for (2.1) proposed in \cite{K1} has been found to be wrong. A counter-example was constructed by Mihara \cite{Mih} who noticed that the matrix
$$
A=\begin{pmatrix}
p & p & 0\\
0 & p & 0\\
0 & 0 & 1\end{pmatrix}
$$
over $\Qp$ or $\Cp$ satisfies the reduction condition from \cite{K1}, though (2.1) is violated for $q(t)=(t-1)(t-p)$.

In fact, the reduction condition from \cite{K1} implies only the equality $\|A^2\|=\|A\|^2$ which is weaker than normality. Only this property is actually proved in \cite{K1} for the operators given there as examples, except the operators related to the non-Archimedean model of the canonical commutation relations of quantum mechanics, whose normality is either immediate from the definition, or follows from a result given below (section 2.3).

\medskip
{\bf 2.2.} {\it Orthoprojections}. A projection on a Banach space $\B$ is such a linear bounded operator $P$ that $P$ is idempotent: $P^2=P$. It is obvious that either $P=0$, or $\|P\|\ge 1$. Kernels of $P\ne 0$ and $I-P$ complement each other having a trivial intersection; if they are orthogonal (in the non-Archimedean sense), then $P$ is called an {\it orthoprojection}. In this and only this case, $\|P\|=1$. For an orthoprojection $P$ different from 0 and $I$, $\|P\|=\|I-P\|$ (\cite{Rooij}, page 63).

\medskip
\begin{teo}
A projection $P$ is an orthoprojection, if and only if it is strongly normal.
\end{teo}

\medskip
{\it Proof.} Note that the spectrum of an orthoprojection is contained in the two-point set $\{ 0,1\}$. Indeed, consider the equation
\begin{equation}
Px-\lambda x=z,\quad z\in \B, \lambda \in K.
\end{equation}
Let us write $z=\xi +\eta$ where $\xi \in M$, $\eta \in N$,
$$
M=\ker (I-P)=\{ x\in \B :\ Px=x\} ,\quad N=\ker P=\{ x\in \B :\ Px=0\}.
$$
We can seek a solution of (2.2) in the form $x=g+h$, $g\in M$, $h\in N$. Then (2.2) means that $g-\lambda g-\lambda h=\xi +\eta$, and if $\lambda \notin \{ 0,1\}$, we can write a solution as follows:
\begin{equation}
g=(1-\lambda )^{-1}\xi ,\quad h=-\lambda^{-1}\eta .
\end{equation}
By (2.3), $P-\lambda I$ is surjective. One verifies directly that $P-\lambda I$ is injective.

Suppose that $P$ is an orthoprojection different from 0 or $I$. In this case $\sigma (P)=\{ 0,1\}$.

The Banach algebra $\LP$ generated by $P$ and $I$ consists of operators $\alpha P+\beta I$, $\alpha ,\beta \in K$, or, equivalently, of operators $aP+b(I-P)$, $a,b\in K$, for which
\begin{equation}
\| aP+b(I-P)\|\le \max (|a|,|b|).
\end{equation}

Let us take such vectors $g\in M$, $h\in N$, that $\|g\|=\|h\|$, and set $x=g+h$. Due to the orthogonality, $\|x\|=\|g\|=\|h\|$. On the other hand, $[aP+b(I-P)]x=ag+bh$, so that
$$
\|[aP+b(I-P)]x\|=\max (|a|,|b|)\cdot \|x\|.
$$
Together with (2.4), this implies the equality
\begin{equation}
\| aP+b(I-P)\|=\max (|a|,|b|).
\end{equation}

It follows from (2.5) that $\LP$ is isomorphic to $C(\{ 0,1\},K)$. Hence $P$ is strongly normal.

Conversely, if $P$ is a strongly normal projection, then under the isomorphism of $\LP$ and $C(\sigma (P),K)$, the operator $P$ corresponds to such a continuous function $f$ that $f^2=f$, so that
$$
f(x)=\begin{cases}
0, & \text{if $x\in X_0$},\\
1, & \text{if $x\in X_1$},
\end{cases}
$$
where $X_0\cup X_1=\sigma (P)$. If $P\ne 0$, then $X_1$ is nonempty, so that $\|P\|=\|f\|=1$, and $P$ is an orthoprojection. \qquad $\blacksquare$

\medskip
{\bf 2.3.} {\it Multiplication operators}. Let us consider the Banach space $\B =C(M,K)$ where $M$ is a compact totally disconnected Hausdorff topological space. It is known that $\B$ possesses an orthonormal basis (\cite{Rooij}, Corollary 5.25, or \cite{PS}, Theorem 2.5.22).

\medskip
\begin{teo}
Let $A$ be an operator of multiplication on $\B$ by a function $a\in \B$. Then $A$ is strongly normal.
\end{teo}

\medskip
{\it Proof}. Obviously, the spectrum $\sigma (A)$ coincides with the set $R=a(M)$ of values of the function $a$. The set $R$ is a compact subset of the zero-dimensional topological field $K$. Therefore (\cite{E}, 6.2.10 and 6.2.11) $R$ is zero-dimensional, hence totally disconnected.

The Banach algebra $\LA$ can be identified with the completion of the set of functions $p(a(m))$, $m\in M$, where $p\in K[r]$, with respect to the norm
$$
\sup\limits_{0\ne \varphi \in \B} \sup\limits_{m\in M}\frac{|p(a(m))\varphi (m)|}{\|\varphi \|}=\sup\limits_{m\in M}|p(a(m))|=\sup\limits_{r\in R}|p(r)|.
$$

By Kaplansky's theorem (Theorem 43.3 in \cite{Sch}), the algebra $\LA$ coincides with the set of operators of multiplication by functions $\pi (a(\cdot))$, $\pi \in C(R,K)$, thus it is isomorphic to $C(R,K)$. $\qquad \blacksquare$

\section{Commutative Algebras and Algebras with Baer Reductions}

{\bf 3.1.} {\it Commutative algebras}. Classically, if a $*$-algebra $\AAA$ of bounded operators on a Hilbert space is commutative, then $AA^*=A^*A$ for any $A\in \AAA$, so that all the operators from $\AAA$ are normal. Therefore in our situation it is reasonable to consider a commutative algebra $\AAA$ of normal operators on a Banach space $\B$ over the field $K$. We assume that $\AAA$ is complete with respect to the norm of operators and contains the unit operator $I$.

\medskip
\begin{teo}
Under the above assumptions, the algebra $\AAA$ is isomorphic to the algebra $C(M,K)$ of $K$-valued continuous functions on a compact totally disconnected Hausdorff topological space $M$. Under this isomorphism, characteristic functions of open-closed subsets $\Lambda \subset M$ correspond to orthoprojections $E(\Lambda )$ forming an orthoprojection-valued finitely additive measure on the algebra of open-closed subsets of $M$. For an operator $F\in \AAA$ corresponding to a function $f\in C(M,K)$, there is an integral representation
\begin{equation}
F=\int\limits_M f(\lambda )\,E(d\lambda )
\end{equation}
convergent with respect to the norm of operators.
\end{teo}

\medskip
{\it Proof}. For any operator $A\in \AAA$, its normality means that the algebra $\LA$ is uniform and all its characters take values in $K$. Then the whole algebra $\AAA$ possesses the same properties; therefore \cite{Ber} $\AAA$ is isomorphic to $C(\mathcal M(\AAA ),K)$. The proof of (3.1) is identical to that of the representation of a single normal operator \cite{K1}. $\qquad \blacksquare$

\medskip
{\bf 3.2.} {\it Algebras with Baer reductions}. In order to introduce a possible non-Archimedean counterpart for the class of von Neumann algebras, we need the following reduction procedure. We assume that the space $\B$ is infinite-dimensional and possesses an orthonormal basis (in the non-Archimedean sense; see \cite{Rob,Rooij}). Then $\B$ is isomorphic to the space $c_0(J,K)$ of sequences $x=(x_1,x_2,\ldots ,x_i,\ldots )$, $i\in J$, $x_i\in K$, $x_i\to 0$, by the filter of complements to finite subsets of an infinite set $J$. We may assume that $
\B =c_0(J,K)$.

Let $\AAA$ be an algebra of linear bounded operators on $c_0(J,K)$. Denote by $\AAA_1$ the closed unit ball in $\AAA$ -- the set of all operators from $\AAA$ with norm $\le 1$. $\AAA_1$ is an algebra over the ring $O$, just as its ideal $\AAA_0$ consisting of operators of norm $<1$. The {\it reduced algebra} $\widetilde{\AAA}=\AAA_1/\AAA_0$ can be considered as a $\widetilde{K}$-algebra. Now we can look for a class of $\widetilde{K}$-algebras, for which there is a (purely algebraic) theory parallel to the theory of von Neumann algebras. Then the class of algebras $\AAA$ corresponding to $\widetilde{\AAA}$ from that class will be the desired one.

An algebraic theory of the above kind is the theory of Baer rings and algebras developed by Kaplansky \cite{Kap}. A unital ring $R$ is called the Baer ring, if each left (or, equivalently, each right) annihilator in $R$ is generated by an idempotent element. This property was proved by Baer \cite{Baer} for the ring of all endomorphisms of a vector space of an arbitrary dimension. Kaplansky \cite{Kap1} proved it for any $AW^*$-algebra (the class of $AW^*$-algebras is wider than the class of von Neumann algebras).

A Baer ring $R$ is called Abelian, if all its idempotents are central, and Dedekind finite, if $xy=1$ implies $yx=1$. An idempotent $e\in R$ is called Abelian (finite), if the Baer ring $eRe$ is Abelian (respectively, Dedekind finite). If $u$ and $v$ are central idempotents, we write $u\le v$, if $vu=u$. An idempotent $e$ is called faithful, if the smallest of the central idempotents $v$ satisfying $ve=e$ is equal to 1.

Kaplansky \cite{Kap} introduced the following types of Baer rings. A Baer ring $R$ is of type I, if it has a faithful Abelian idempotent. It is of type II, if it has a faithful finite idempotent, but no nonzero Abelian idempotents, and of type III, if it has no nonzero finite idempotents. These classes are subdivided further into finite and infinite ones.

A typical example of a type I Baer ring is the ring of all linear transformations of a vector space of countable dimension. Under some additional conditions, an arbitrary Baer factor of type I is of this form \cite{Wo2}. For other examples of Baer rings see \cite{Kap,Ho1,Ho2,Wo1}.

The main result \cite{Kap} of the theory of Baer rings is the unique decomposition of every Baer ring into a direct sum of rings of the above types.

Let us call an operator algebra $\AAA$ {\it an algebra with the Baer reduction}, if the reduced algebra $\AR$ is a Baer ring. It is well known (see \cite{BGR}, Lemma 2.5.1/3) that a finite system of elements of norm 1 in a non-Archimedean normed space is orthonormal if and only if their reductions are linearly independent. Therefore {\bf the operator ring $\AAA_1$ with the Baer reduction is an orthogonal sum of rings with reductions of types I, II, and III}.

The simplest example of algebras with the Baer reduction of type I is the algebra of all bounded operators on a Banach space of countable type (that is, $c_0(J,K)$ with a countable set $J$). The very existence of other operator algebras (moreover, factors) with this property is far from obvious. Below we present a class of such algebras. As a whole, the study of various classes of operator algebras with the Baer reduction seems a huge problem comparable with the whole theory of von Neumann algebras.

\section{Non-Archimedean Crossed Products}

{\bf 4.1.} {\it Analysis on product spaces}. Let $S$ be a totally disconnected compact Hausdorff topological space, $G$ be an Abelian infinite second-countable totally disconnected compact Hausdorff topological group acting transitively on $S$ by homeomorphisms. The action will be denoted as $x\mapsto xa$, $x\in S$, $a\in G$ (the group operation is written multiplicatively). Below we construct and study some operators on $C(S\times G,\Cp)$, the space of continuous functions on $S\times G$ with values in $\Cp$. Here $p$ is a prime number, $\Cp$ is the completion of an algebraic closure of the field $\Qp$ of $p$-adic numbers. We will denote $|\cdot |_p$ the absolute value in $\Cp$. Some features of our approach follow \cite{MN} though the actual meaning of our objects is different. Note also (in order to avoid confusion referring to various sources) that in the class of compact Hausdorff spaces the properties of total disconnectedness and zero dimensionality are equivalent.

We begin with the following easy property of the space $C(S,\Cp )$. As usual, we denote by $\AAA'$ the commutant of an operator algebra $\AAA$ over a Banach space $\B$, that is the set of such bounded operators $B$ on $\B$ that $AB=BA$ for each $A\in \AAA$. Below we use a similar notation for the commutant of an arbitrary set of operators.

\medskip
\begin{lem}
Let $L_\varphi$ be the operator of multiplication by a function $\varphi \in C(S,\Cp)$ acting on $C(S,\Cp)$, $\mathbf L$ be the algebra of all such operators. Then $\mathbf L'=\mathbf L$.
\end{lem}

\medskip
{\it Proof}. Obviously, $\mathbf L\subset \mathbf L'$. Let $A\in \mathbf L'$. Then, for every $\varphi \in C(S,\Cp)$, $AL_\varphi =L_\varphi A$, that is $A(\varphi f)=\varphi (Af)$ for any $f\in C(S,\Cp)$. For $f=1$, we get $(A\varphi )(x)=\varphi (x)(A1)(x)$, so that $A=L_\psi$, $\psi (x)=(A1)(x)$, $x\in S$. Therefore $A\in \mathbf L.\qquad \blacksquare$

\medskip
From now on, we assume that {\bf the group $G$ is $p$-compatible}. This notion is defined as follows (see \cite{Rooij}). Denote by $o(G)$ the set of all such natural numbers $n$, for which there exists a subgroup $H\subset G$ with the property that $G/H$ has an element of order $n$. The $p$-compatibility means that $p\notin o(G)$ (an additional condition formulated in \cite{Rooij}, Sect. 9.J, for more general fields is satisfied automatically for $\Cp$). A typical example is $G=S=\mathbb Z_l$ where $l$ is a prime different from $p$.

It is well known (see \cite{Rooij,Sch71}) that $p$-compatible groups possess nice properties resembling those appearing in classical harmonic analysis. In particular, $G$ possesses a $\Cp$-valued Haar measure $\mu$. In the Banach space $C(G,\Cp )$, there is an orthonormal basis $\{ g_j\}$ consisting of $\Cp$-valued characters, that is maps $G\to \Cp$, such that $g_j(ab)=g_j(a)g_j(b)$. It is convenient to index the characters not by natural numbers but by elements of the dual group $\GD$ consisting of all $\Cp$-valued characters or, equivalently, of all continuous homomorphisms $G\to \mathbb T_p$ where $\mathbb T_p$ is the set of all roots of 1 in $\Cp$ of orders prime to $p$ (\cite{Rooij}, p. 360). For groups of this type, $\GD$ is isomorphic also to the Pontryagin dual (\cite{Rooij}, p. 350). The second-countability property of $G$ implies its metrizability and the countability of $\GD$ (\cite{HR}, Theorem 24.15). Below we will use the fact that $\GD$ (with the discrete topology) is torsional, that is every finite subset of $G$ lies in a finite subgroup (see Corollary 9.15 in \cite{Rooij}).

Returning to characters $g_j$ we stress that $g_i$, $i\in \GD$, is identical with $i$, so that $g_i(a)=i(a)$, for any $a\in G$. In particular, $g_1(a)\equiv 1$. Thus, $g_i(a)g_j(a)=g_{i\cdot j}(a)$, $g_i(a^{-1})=g_{i^{-1}}(a)$. Note also \cite{Sch71} that the characters are orthonormal not only in the non-Archimedean sense, but also in the integral sense, with respect to the $\Cp$-valued Haar measure $\mu$:
\begin{equation}
\int\limits_G g_j(a)g_n(a^{-1})\,\mu (da)=\delta_{j,n},\quad j,n\in \GD,
\end{equation}
where $\delta_{j,n}$ is the Kronecker symbol.

Dealing with a function $F\in C(S\times G,\Cp )$ we can write
\begin{equation}
F(x,a)=\sum\limits_{n\in \GD }\varphi_n(x)g_n(a),\quad x\in S,a\in G,
\end{equation}
where
\begin{equation}
\varphi_n(x)=\int\limits_G F(x,a)g_n(a^{-1})\,\mu (da)
\end{equation}
(for the integration theory with non-Archimedean-valued measures see \cite{Rooij,Sch}). It follows from (4.3) that the functions $\varphi_n$ are continuous. Below we use the notation $F\sim \langle \varphi_n(x)\rangle _{n\in \GD}$ where $\|\varphi_n\|\to 0$ by the filter of complements to finite sets in $\GD$ (in such cases we will write $n\to \infty$); see, for example, Proposition 1.6 in \cite{K09}. Here and below we denote the supremum norms in various spaces of continuous functions by the same symbol $\|\cdot \|$.

Let $\overline{A}$ be a linear bounded operator on $C(S\times G,\Cp )$ (letters with a bar will denote operators on $C(S\times G,\Cp )$, to distinguish them from operators on $C(S,\Cp )$). Then for any $\varphi \in C(S,\Cp )$, $n\in \GD$,
\begin{equation}
(\overline{A}(\varphi g_n))(x,a)=\sum\limits_{j\in \GD }y_{j,n}(x)g_j(a)
\end{equation}
where $y_{j,n}\in C(S,\Cp )$, $\|y_{j,n}\| \to 0$, as $j\to \infty$,
\begin{equation}
\sup\limits_{j\in \GD}\|y_{j,n}\|=\|\overline{A}(\varphi g_n)\|\le \|\overline{A}\| \cdot \|\varphi \|.
\end{equation}

Define an operator $A_{j,n}$ on $C(S,\Cp )$ by setting $A_{j,n}\varphi =y_{j,n}$. By (4.4) and (4.5), this operator is linear and bounded. For each $n$, $A_{j,n}\to 0$ in the strong operator topology. By (4.5),
\begin{equation}
\sup\limits_{j,n\in \GD}\|A_{j,n}\|\le \|\overline{A}\|.
\end{equation}

It follows from (4.2) and (4.4) that
\begin{equation}
\overline{A}F=\sum\limits_{n\in \GD }\overline{A}(\varphi_ng_n)=\sum\limits_n\sum\limits_j(A_{j,n}\varphi_n)g_j
=\sum\limits_j g_j\sum\limits_nA_{j,n}\varphi_n.
\end{equation}
The change in the order of summation is justified as follows. By Theorem 3.8 from \cite{Ka} (which remains valid in the vector case and for the summation over $\GD$), it is sufficient to prove that, for any $\varepsilon >0$, there exists a finite subset $N=N(\varepsilon )\subset \GD$, such that $\|(A_{j,n}\varphi_n)g_j\| <\varepsilon$ whenever $n,j\notin N$.

Indeed, we know that $\|\varphi_n\| \to 0$. By (4.6), $\|A_{j,n}\varphi_n\|\le \|\overline{A}\| \cdot \|\varphi_n\|$.  Choose a finite set $M$ in such a way that $\|\varphi_n\|<\frac{\varepsilon}{\|\overline{A}\|}$, for $n\notin M$. On the other hand, for a fixed $n$, $A_{j,n}\varphi_n$ is the $j$-th coefficient of the expansion of the function $\overline{A}(\varphi_ng_n)$ on $G$ in the basis $\{g_j\}$. Therefore we find such a finite subset $J\subset \GD$ that
$$
\|A_{j,n}\varphi_n\|<\varepsilon ,\quad n\in M,
$$
if $j\notin J$. Thus, $\|A_{j,n}\varphi_n\|<\varepsilon$ for all $n\in \GD$, if $j\notin J$, so that
$$
\|(A_{j,n}\varphi_n)g_j\| <\varepsilon \quad \text{for all $n\in \GD,j\notin J$}.
$$
Meanwhile,
$$
\|(A_{j,n}\varphi_n)g_j\| <\varepsilon \quad \text{for all $n\notin M,j\in \GD$}.
$$
Therefore we obtain the required inequality, if we set $N=M\cup J$.

The equalities (4.7) show that the correspondences $F\sim \langle \varphi_n\rangle$, $\overline{A}\sim \langle A_{j,n}\rangle$, agree with the usual rules of linear algebra. Similarly, sums and products of these infinite matrices correspond to the sums and products of operators.

\bigskip
{\bf 4.2.} {\it Approximation on $\GD$}. Below we will need a certain approximation property for functions on $\GD$ with values in $\Cp$.

Let $\gamma :\ \GD \to \mathbb R_+$ be a fixed function tending to zero at infinity. Denote by $l_\gamma$ the Banach space of functions $f:\ \GD\to \Cp$ with finite norm
$$
\|f\|_\gamma =\sup\limits_{i\in \GD}|f(i)|_p\gamma (i).
$$
We call a function $i\mapsto \sum\limits_a c_ag_i(a)$, where $c_a\in \Cp$ and the sum is indexed by a finite set of elements $a\in G$, {\it a trigonometric polynomial} on $\GD$.

\medskip
\begin{lem}
Any bounded function $f:\ \GD\to \Cp$ can be arbitrarily well approximated in $l_\gamma$ by trigonometric polynomials.
\end{lem}

\medskip
{\it Proof}. Suppose that $\sup\limits_{i\in \GD}|f(i)|_p=M_f<\infty$. Given $\varepsilon >0$, choose a finite subgroup $\Sigma \subset \GD$ in such a way that
\begin{equation}
\sup\limits_{i\notin \Sigma}|\gamma (i)|<\frac{\varepsilon}{M_f}
\end{equation}
(that is possible, since the group $\GD$ is torsional). Let us consider the restriction of $f$ to the subgroup $\Sigma$. Since $\Sigma$ is finite, this restriction can be expanded into a finite linear combination of characters on $\Sigma$ with values in $\mathbb T_p$. This is a consequence of the fact that $\GD$ and all its subgroups are $p$-compatible (\cite{Rooij}, p. 361). Note also that for such groups there is an isomorphism between $\mathbb C$-valued and $\mathbb T_p$-valued characters (\cite{Rooij}, p. 360). In addition, every character of $\Sigma$ can be extended to a character on $\GD$ (\cite{Rooij}, p. 350). Since $(\GD )\sphat =G$, each character on $\GD$ has the form $i\mapsto g_i(a)$ where $a\in G$ is a fixed element. Thus, we have the representation
$$
f(i)=\sum\limits_{a\in \widehat{\Sigma}}c_ag_i(a),\quad i\in \Sigma ,
$$
where $c_a\in \Cp$,
$$
\sup\limits_{i\in \Sigma}|f(i)|_p=\sup\limits_{a\in \widehat{\Sigma}}|c_a|_p,
$$
so that $|c_a|_p\le M_f$ for $a\in \widehat{\Sigma}$.

Set
$$
f_\varepsilon (i)=\sum\limits_{a\in \widehat{\Sigma}}c_ag_i(a),\quad i\in \GD .
$$
Then $f_\varepsilon (i)-f(i)=0$ for $i\in \Sigma$,
$$
\|f_\varepsilon -f\|_\gamma =\sup\limits_{i\ne \Sigma}|f_\varepsilon (i)-f(i)|_p\gamma (i)<\varepsilon
$$
by virtue of (4.8), since $|g_i(a)|_p\le 1$, so that $|f_\varepsilon (i)|_p\le M_f.\qquad \blacksquare$

\medskip
Note that the dual space to the Banach space $l^\infty (\GD ,\Cp )$ of bounded $\Cp$-valued functions on $\GD$ is $c_0(\GD ,\Cp )$ (see Theorem 5.5.5 in \cite{PS}). Therefore the above reasoning proves also that the set of all trigonometric polynomials is weakly dense in $l^\infty (\GD ,\Cp )$.

\bigskip
{\bf 4.3.} {\it The crossed product construction}. On the Banach space $\B=C(S\times G,\Cp )$, we consider the operators
\begin{gather*}
\overline{U}_{a_0}F(x,a)=F(xa_0,aa_0);\\
\overline{V}_{a_0}F(x,a)=F(x,a_0^{-1}a);\\
\overline{W}F(x,a)=F(xa^{-1},a^{-1});\\
\overline{L}_\varphi F(x,a)=\varphi (x)F(x,a);\\
\overline{M}_\varphi F(x,a)=\varphi (xa^{-1})F(x,a),
\end{gather*}
$x\in S$, $a\in G$. Here $a_0\in G$ is a fixed element, $\varphi \in C(S,\Cp )$ is a fixed function.

All the above operators are bounded. It is easy to check that
\begin{equation}
\overline{W}=\overline{W}^{-1},\ \overline{W}\, \overline{U}_{a_0}\overline{W}=\overline{V}_{a_0},\ \overline{W}\, \overline{L}_\varphi \overline{W}=\overline{M}_\varphi .
\end{equation}

Denote by $\II$ the set of all the operators $\overline{U}_{a_0}$ and $\overline{L}_\varphi$, and by $\JJ$ the set of all the operators $\overline{V}_{a_0}$ and $\overline{M}_\varphi$. The closed linear hull (with respect to the strong operator topology \cite{KR,Ta}) of the set $\II$ will be denoted $\RR (\II )$. Similarly, $\RR (\JJ )$ is the strongly closed linear hull of $\JJ$. It follows from (4.9) that the mapping $\overline{A}\mapsto \overline{W}\, \overline{A}\, \overline{W}$ is a spatial isomorphism of the algebras $\RR (\II )$ and $\RR (\JJ )$. These algebras will be interpreted as non-Archimedean crossed product algebras, at least for the case of $p$-compatible groups.

Let us compute matrix representations of the above operators. Let $F$ be given by (4.2). Then
$$
\left( \overline{U}_{a_0}F\right) (x,a)=\sum\limits_{n\in \GD}g_n(a_0)\left( U_{a_0}\varphi_n\right) (x)g_n(a)
$$
where $U_{a_0}$ is an operator on $C(S,G)$, $\left( U_{a_0}f\right) (x)=f(xa_0)$. Therefore
\begin{equation}
\overline{U}_{a_0}\sim \langle \delta_{j,n}g_n(a_0)U_{a_0}\rangle _{j,n\in \GD}.
\end{equation}

If $\overline{A}\sim \langle A_{j,n}\rangle$ is an arbitrary bounded operator, then by (4.10),
\begin{equation}
\overline{U}_{a_0}^{-1}\overline{A}\, \overline{U}_{a_0}\sim \langle g_j(a_0^{-1})g_n(a_0)U_{a_0}^{-1}
A_{j,n}U_{a_0}\rangle.
\end{equation}
Also we have $\overline{L}_\varphi \overline{A}\sim \langle L_\varphi A_{j,n}\rangle$, $\overline{A}\, \overline{L}_\varphi \sim \langle A_{j,n}L_\varphi \rangle$. Thus, if an operator $\overline{A}$ commutes with all $\overline{L}_\varphi$, then by Lemma 4.1, each operator $A_{j,n}$ is the operator of multiplication by a continuous function $\psi_{j,n}$ on $S$. If $\overline{A}\in \II'$, then, in addition, $\overline{A}=\overline{U}_{a_0}^{-1}\overline{A}\, \overline{U}_{a_0}$, and we find from (4.11) that for any $f\in C(S,\Cp )$,
$$
\psi_{j,n}(x)f(x)=g_j(a_0^{-1})g_n(a_0)\psi_{j,n}(xa_0^{-1})f(x),
$$
which implies the identity
\begin{equation}
\psi_{j,n}(xa_0)=g_j(a_0^{-1})g_n(a_0)\psi_{j,n}(x),\quad \text{for any $x\in S,a_0\in G$}.
\end{equation}

Let us find a representation of operators from $\JJ$. We have, for any $\varphi \in C(S,\Cp )$,
$$
\overline{V}_{a_0}(\varphi (x)g_n(a))=\varphi (x)g_n(a_0^{-1})g_n(a),
$$
so that
\begin{equation}
\overline{V}_{a_0}\sim \langle \delta_{j,n}g_j(a_0^{-1})\rangle .
\end{equation}
Next,
$$
\overline{M}_\psi (\varphi (x)g_n(a))=\psi (xa^{-1})\varphi (x)g_n(a),\quad \psi \in C(S,\Cp ).
$$
Writing
\begin{equation}
\psi (xa^{-1})=\sum\limits_{m\in \GD}c_m(x)g_m(a)
\end{equation}
we get
$$
\overline{M}_\psi (\varphi (x)g_n(a))=\sum\limits_m c_m(x)\varphi (x)g_{n\cdot m}(a),
$$
so that
\begin{equation}
\overline{M}_\psi \sim \langle M_{j,n}\rangle
\end{equation}
where $M_{j,n}$ is the operator of multiplication by the function $c_{n^{-1}\cdot j}$, an appropriate coefficient from (4.14).

Let us fix an element $x_0\in S$. Denote by $\GD_0$ the set of such $i\in \GD$ that the equality $x_0a=x_0$ implies the equality $g_i(a)=1$. The set $\GD_0$ is a subgroup in $\GD$. If the action is free, that is $x_0a=x_0$ implies $a=1$, then $\GD_0=\GD$.

For $i\in \GD_0$, define a function $\eta_i\in C(S,\Cp )$ as follows. Since the action of $G$ on $S$ is transitive, for each $x\in S$ there exists an element $a\in G$, such that $x=x_0a^{-1}$. Set $\eta_i(x)=g_i(a)$. If also $x=x_0b^{-1}$, $b\in G$, then $x_0a^{-1}b=x_0$, and $g_i(a)=g_i(b)$, provided $i\in \GD_0$. Therefore the function $\eta_i$ is well-defined. Note also that $\|\eta_i\|=1$ and $\eta_1(x)\equiv 1$.

As we have proved, an operator $\overline{A}\in \II'$ is represented by the matrix of multiplication operators corresponding to the functions $\psi_{j,n}$ on $S$, satisfying (4.12) or, equivalently, the identity
\begin{equation}
\psi_{j,n}(xa^{-1})=g_{j\cdot n^{-1}}(a)\psi_{j,n}(x),\quad x\in S,a\in G.
\end{equation}
Set $x=x_0$. If $i=j\cdot n^{-1}\notin \GD_0$, then there exists such an element $a\in G$ that $x_0a^{-1}=x_0$ and $g_i(a)\ne 1$. It follows from (4.16) that in this case $\psi_{j,n}(x_0)=0$, so that $\psi_{j,n}(x)=0$ for all $x\in S$. Therefore
\begin{equation}
\psi_{j,n}(x)=\begin{cases}
b_{j,n}\eta_{j\cdot n^{-1}}(x), & \text{if $j\cdot n^{-1}\in \GD_0$};\\
0, & \text{if $j\cdot n^{-1}\notin \GD_0$},
\end{cases}
\end{equation}
where $b_{j,n}=\psi_{j,n}(x_0)$.

The function $\eta_i$ also enables us to give a more explicit representation of the coefficients $c_{j\cdot n^{-1}}(x)$ appearing in (4.15). We find, using the orthogonality relation (4.1), that
$$
c_{j\cdot n^{-1}}(x)=\int\limits_G\psi (xa^{-1})g_{j\cdot n^{-1}}(a^{-1})\,\mu (da).
$$
For any $\gamma \in G$, using the invariance of the Haar measure, we get
$$
c_{j\cdot n^{-1}}(x\gamma^{-1})=\int\limits_G\psi (xb^{-1})g_{j\cdot n^{-1}}(\gamma b^{-1})\,\mu (db)=g_{j\cdot n^{-1}}(\gamma )c_{j\cdot n^{-1}}(x),
$$
and setting $x=x_0$ we find that
\begin{equation}
c_{j\cdot n^{-1}}(x)=c_{j\cdot n^{-1}}(x_0)\eta_{j\cdot n^{-1}}(x),\quad x\in S,
\end{equation}
where $c_{j\cdot n^{-1}}(x_0)=0$, if $j\cdot n^{-1}\notin \GD_0$.

In particular, if $\psi =\eta_l$, $l\in \GD_0$, then the representation (4.14) takes the form (with $b\in G$ substituted for $a$)
$$
\eta_l(xb^{-1})=\eta_l(x_0a^{-1}b^{-1})=g_l(ab)=\eta_l(x)g_l(b),
$$
so that in the representation (4.14) for $\psi =\eta_l$,
\begin{equation}
c_m(x)=\delta_{l,m}\eta_l(x).
\end{equation}

\bigskip
{\bf 4.4.} {\it Commutants}. We use the above matrix representations to find $\II'=\RR (\II )'$ and $\JJ'=\RR (\JJ )'$.

The next result shows that the algebras $\RR (\II )$ and $\RR (\JJ )$ possess properties resembling those of von Neumann algebras.

\medskip
\begin{teo}
\begin{description}
\item[(i)] $\RR (\JJ )=\II'$; $\RR (\II )=\JJ'$;

\item[(ii)] $\RR (\JJ )=\RR (\JJ )''$; $\RR (\II )=\RR (\II )''$;

\item[(iii)] If the action of $G$ on $S$ is free, then $\RR (\II )$ and $\RR (\JJ )$ are factors.
\end{description}
\end{teo}

\medskip
{\it Proof}. Due to the spatial isomorphism (4.9), it is sufficient to prove one property in each pair.

It follows directly from the definitions that $\JJ \subset \II'$. Since $\II'$ is closed in $L(\B )$ with respect to the strong operator topology, the equality $\RR (\JJ )=\II'$ will be proved if we show that each element from $\II'$ can be approximated in strong operator topology by linear combinations of operators from $\JJ$.

As we know, every operator $\overline{A}$ from $\II'$ has a matrix representation $\langle \psi_{j\cdot n^{-1}}\rangle_{j,n\in \GD}$ where $\psi_{j\cdot n^{-1}}$ is the operator of multiplication by the function (4.17), in which $|b_{j,n}|_p\le C$ for all $j,n$ (here and below the letter $C$ denotes various positive constants), $b_{j,n}\to 0$, as $j\to \infty$.

Let $l\in \GD$, and $\overline{A}_l$ be an operator with the same matrix elements as $\overline{A}$ on a single ``diagonal'' $\{ (j,n):\ j\cdot n^{-1}=l\}$ and zeroes on the rest of matrix entries. Let us prove that $\overline{A}_l\to 0$ in the strong operator topology, as $l\to \infty$.

Let $F\in \B$, $F\sim \langle \varphi_n(x)\rangle$. Then
$$
\overline{A}_lF\sim \left\langle \sum\limits_{n:\ j\cdot n^{-1}=l}b_{j,n}\eta_{j\cdot n^{-1}}(x)\varphi_n(x)\right\rangle_{j\in \GD}=\left\langle b_{j,j\cdot l^{-1}}\eta_l\varphi_{j\cdot l^{-1}}\right\rangle_{j\in \GD},
$$
if $l\in \GD_0$, and $\overline{A}_lF=0$, if $l\notin \GD_0$. It is sufficient to consider the case where $l\in \GD_0$. For any $\varepsilon >0$, there exists such a finite set $\Sigma \subset \GD$ that $\|\varphi_n\|< \dfrac{\varepsilon}C$, as $n\in \GD \setminus \Sigma$. If $j\notin l\cdot \Sigma$, then $\|\varphi_{j\cdot l^{-1}}\|< \dfrac{\varepsilon}C$ and $\|b_{j,j\cdot l^{-1}}\eta_l\varphi_{j\cdot l^{-1}}\|<\varepsilon$. If $j\in l\cdot \Sigma$, then $|b_{j,j\cdot l^{-1}}|_p< \dfrac{\varepsilon}C$, as $j$ is outside a finite set $P=\{ j_1,\ldots ,j_k\}$ (because $j\cdot l^{-1}$ belongs to the finite set $\Sigma$).

Let $Q=\{ m^{-1}j_\nu :\ m\in \Sigma ,\nu =1,\ldots ,k\}$. If $l\notin Q$ and $j\in l\cdot \Sigma$, then $j\notin P$. Therefore, if $l\notin Q$, then
$$
\|b_{j,j\cdot l^{-1}}\eta_l\varphi_{j\cdot l^{-1}}\|<\varepsilon ,
$$
so that $\overline{A}_l\to 0$ in the strong operator topology.

The ultrametric property shows that the operator $\overline{A}$ can be approximated in the strong operator topology by finite sums of operators $\overline{A}_l$. It remains to prove that each $\overline{A}_l$ can be approximated in the same topology by linear combinations of operators from $\JJ$. Thus, from now on, $l$ will be fixed.

Given $\varepsilon >0$ and a finite collection of functions $F_1,\ldots ,F_k\in \B$, we have to find such an operator $\overline{B}_l\in \JJ$ that $\left\|(\overline{A}_l-\overline{B}_l)F_i\right\| <\varepsilon$, $i=1,\ldots ,k$.

Let $F_i\sim \left\langle \varphi_n^{(i)}\right\rangle_{n\in \GD}$. Set
$$
\gamma (j)=\max\limits_{1\le i\le k}\|\varphi^{(i)}_{j\cdot l^{-1}}\| ,\quad j\in \GD .
$$
By Lemma 4.2, there exists such a trigonometric polynomial $i\mapsto \sum\limits_{\nu =1}^N  d_\nu g_i(a_\nu^{-1})$, $d_\nu \in \Cp$, $a_\nu \in G$, that
\begin{equation}
\sup\limits_{i\in \GD}\left| b_{i,i\cdot l^{-1}}-\sum\limits_{\nu =1}^N d_\nu g_i(a_\nu^{-1})\right|_p\gamma (i)<\varepsilon .
\end{equation}

Set
$$
\overline{B}_l=\sum\limits_{\nu =1}^N d_\nu \overline{M}_{\eta_l}\overline{V}_{a_\nu}.
$$
Identifying an operator of multiplication by a certain function with that function and using (4.13), (4.15), and (4.19), we get for the ``matrix elements'' the expression
$$
\left( \overline{A}_l-\sum\limits_{\nu =1}^N d_\nu \overline{M}_{\eta_l}\overline{V}_{a_\nu}\right)_{j,n}=\delta_{l,j\cdot n^{-1}}\left[  b_{j,j\cdot l^{-1}}-\sum\limits_{\nu =1}^N d_\nu g_i(a_\nu^{-1})\right] \eta_l(x),
$$
whence
$$
(\overline{A}_l-\overline{B}_l)F_i\sim \left\langle \left[ b_{j,j\cdot l^{-1}}-\sum\limits_{\nu =1}^N d_\nu g_i(a_\nu^{-1})\right]\varphi_{j\cdot l^{-1}}^{(i)}(x)\eta_l(x)\right\rangle_{j\in \GD},\quad i=1,\ldots ,k,
$$
and it follows from (4.20) that $\left\|(\overline{A}_l-\overline{B}_l)F_i\right\| <\varepsilon$. This proves (i).

It follows from the first equality in (i) that $\II''=\RR (\JJ )'$. On the other hand, $\RR (\JJ )'=\JJ'=\RR (\II )$, by the second equality in (i). Therefore $\RR (\II )=\II''=\RR (\II )''$. Similarly, $\RR (\JJ )=\RR (\JJ )''$.

Let us consider the case of free action and prove that $\RR (\II )$ is a factor. Since $\RR (\II )=\JJ'$, it is sufficient to prove that $\II'\cap \JJ'=\{\lambda I,\lambda \in \Cp\}$.

Suppose that an operator $\overline{A}\sim \langle \psi_{j,n}(x)\rangle$ from $\II'$ commutes with $\JJ$, that is with operators $\overline{M}_\varphi$ for all $\varphi \in C(S,\Cp )$, and with operators $\overline{V}_{a_0}$ for all $a_0\in G$. By (4.13),
$$
 \left( \overline{A}\overline{V}_{a_0}\right)_{i,j}=\psi_{i,j}(x)g_j(a_0^{-1}),
 $$
$$
  \left( \overline{V}_{a_0}\overline{A}\right)_{i,j}=g_i(a_0^{-1})\psi_{i,j}(x).
 $$
Therefore the commutation of $\overline{A}$ and $\overline{V}_{a_0}$ for all $a_0\in G$ means that $\psi_{i,j}(x)\equiv 0$ if $i\ne j$. On the other hand, it follows from (4.12) and the transitivity of the action of $G$ on $S$ that $\psi_{i,i}(x)\equiv b_i\in \Cp$, for each $i\in \GD$.

Now we consider the commutation with $\overline{M}_\varphi$ using (4.15) and (4.19). We have $\GD_0=\GD$ by our assumption, so that
$$
 \left( \overline{A}\overline{M}_\varphi \right)_{i,j}=b_i\eta_{i\cdot j^{-1}}(x);
 $$
$$
 \left( \overline{M}_\varphi \overline{A}\right)_{i,j}=b_j\eta_{i\cdot j^{-1}}(x),
 $$
and we get $b_i=b_j$, which means that $\overline{A}=\lambda I.\qquad \blacksquare$

\medskip
Note that in the general situation of a non-free action of $G$, $\RR (\II )$ and $\RR (\JJ )$ are ``almost'' factors -- their central elements have the matrix representation $\langle \delta_{j,n}\lambda_jI\rangle$ where $\lambda_j$ are constants whose values coincide for $j\in \GD_0$.

\medskip
\begin{cor}
Idempotents in $\RR (\JJ )$ have the form $\overline{A}\sim \langle \psi_{j,n}(x)\rangle$ where
$$
\psi_{j,n}(x)=\begin{cases}
b_{j,n}\eta_{j\cdot n^{-1}}(x), & \text{if $j\cdot n^{-1}\in \GD_0$},\\
0, & \text{if $j\cdot n^{-1}\notin \GD_0$},\end{cases}
$$
$b_{j,n}$ are elements from $\Cp$, such that
$$
\sum\limits_{n:\substack{i\cdot n^{-1}\in \GD_0\\n\cdot j^{-1}\in \GD_0}}b_{i,n}b_{n,j}=b_{i,j},\quad i,j\in \GD_0.
$$
The above idempotents are orthoprojections if and only if $|b_{j,n}|_p\le 1$, for all $j,n\in \GD_0$.
\end{cor}

\medskip
{\it Proof}. We obtain the required matrix representation from the fact that $\RR (\JJ )=\II'$ using the equality (4.17). $\qquad \blacksquare$

\bigskip
{\bf 4.5.} {\it The Baer property}. In order to study the reduction of operators from the unit ball $\AAA_1$ of the algebra $\AAA =\RR (\JJ )$, we will need an orthonormal basis in $C(S\times G,\Cp )$.

\medskip
\begin{lem}
The collection of products $\eta_i(x)g_j(a)$ ($x\in S$, $a\in G$), $i\in \GD_0$, $j\in \GD$, forms an orthonormal basis in $C(S\times G,\Cp )$.
\end{lem}

\medskip
{\it Proof}. It is sufficient to show that $\{ \eta_i\}_{i\in \GD_0}$ is an orthonormal basis in $C(S,\Cp )$ (see \cite{Amice}).

The mapping $C(S,\Cp )\to C(G,\Cp )$, $f(x)\mapsto f(x_0a^{-1})$, is an isometric imbedding. Its image $X$ consists of those functions $\varphi$, for which $\varphi (a)=\varphi (b)$ if $x_0a^{-1}=x_0b^{-1}$ (given such a function $\varphi$, one can restore $f$ setting $f(x)=\varphi (a)$ where $x=x_0a^{-1}$).

Let $\varphi \in X$. Then in the expansion
$$
\varphi =\sum\limits_{i\in \GD}c_ig_i,\quad c_i\in \Cp,
$$
only the coefficients $c_i$, $i\in \GD_0$, can be different from 0. Indeed, if $i\notin \GD_0$, then there exists $b\in G$, such than $x_0b^{-1}=x_0$, $g_i(b)\ne 1$. Then $x_0b^{-1}a^{-1}=x_0a^{-1}$, $\varphi (ba)=\varphi (a)$, $a\in G$. For this $i$, due to the invariance of the Haar measure,
\begin{multline*}
c_i=\int\limits_G\varphi (a)g_i(a^{-1})\,\mu (da)=\int\limits_G\varphi (ba)g_i(b^{-1}a^{-1})\,\mu (da)\\
 =[g_i(b)]^{-1}\int\limits_G\varphi (a)g_i(a^{-1})\,\mu (da)=[g_i(b)]^{-1}c_i,
\end{multline*}
so that $c_i=0$.

Thus, the characters $g_i$, $i\in \GD_0$, form an orthonormal basis in $X$. It follows that $\{ \eta_i\}_{i\in \GD_0}$ is an orthonormal basis in $C(S,\Cp )$. $\qquad \blacksquare$

\medskip
Let $\AR$ be the reduced algebra (over the algebraic closure $\mathcal F$ of the finite field $\mathbb F_p$) corresponding to the algebra $\AAA =\RR (\JJ )$; see Section 3.2.

\medskip
\begin{teo}
$\AR$ is a type I Baer ring.
\end{teo}

\medskip
{\it Proof}. Let us use the orthonormal basis in $C(S\times G,\Cp )$ given in Lemma 4.3. The basis elements $\nu_{(i,j)}=\eta_i(x)g_j(a)$ are indexed by pairs $(i,j)$, $i\in \GD_0$, $j\in \GD$. If $\overline{A}\in \RR (\JJ )=\II'$, then by (4.17),
\begin{equation}
\overline{A}(\eta_i(x)g_j(a))=\sum\limits_{\substack{m\in \GD \\ m\cdot j^{-1}\in \GD_0}}b_{m,j}\eta_{m\cdot i\cdot j^{-1}}(x)g_m(a)
\end{equation}
where $b_{m,j}\in \Cp$. Writing $b_{m,j}=0$ for $m\cdot j^{-1}\notin \GD_0$ we can rewrite (4.21) as the matrix representation
$$
\overline{A}\nu_{(i,j)}=\sum\limits_{(l,m)\in \GD_0\times \GD}c_{(i,j),(l,m)}\nu_{(l,m)}
$$
where
\begin{equation}
c_{(i,j),(l,m)}=b_{m,j}\delta_{l,m\cdot i\cdot j^{-1}}.
\end{equation}
The operators $\overline{A}$ from the unit ball $\AAA_1$ are those with $|b_{m,j}|_p\le 1$.

It follows from (4.22) that the matrix representation of a reduced operator from $\AR$ has matrix elements
\begin{equation}
\widetilde{c}_{(i,j),(l,m)}=\widetilde{b}_{m,j}\delta_{l,m\cdot i\cdot j^{-1}},\quad (i,j),(l,m)\in \GD_0\times \GD,
\end{equation}
where $\widetilde{b}_{m,j}\in \mathcal F$, $\widetilde{b}_{m,j}=0$ if $m\cdot j^{-1}\notin \GD_0$. It is easy to check that a product of two such matrices corresponds to the product of the matrices given by the first factors on the right in (4.23).

Hence, the algebra $\AR$ is isomorphic to the algebra $\widetilde{\mathfrak B}$ of column-finite matrices $\left( \widetilde{b}_{m,j}\right)_{m,j\in \GD}$ over $\mathcal F$, such that $\widetilde{b}_{m,j}=0$ if $m\cdot j^{-1}\notin \GD_0$. Such matrices can be interpreted as operators $\widetilde{B}$ on the $\mathcal F$-vector space $Z$ of sequences $(z_i)_{i\in \GD}$ containing only a finite number of nonzero elements, with the following additional property. Let $Z_0\subset Z$ consist of such sequences that $z_i=0$ for $i\notin \GD_0$, while $Z_1$ is the set of such sequences that $z_i=0$ for $i\in \GD_0$. Then the subspaces $Z_0$ and $Z_1$ are invariant for the operators $\widetilde{B}$.

Indeed, if $z=(z_i)\in Z_0$, then
$$
(\widetilde{B}z)_m=\sum\limits_{i\in \GD_0}\widetilde{b}_{m,j}z_i=0,\quad \text{if $m\notin \GD_0$}.
$$
If $z\in Z_1$, then
$$
(\widetilde{B}z)_m=\sum\limits_{i\notin \GD_0}\widetilde{b}_{m,j}z_i=0,\quad \text{if $m\in \GD_0$}.
$$

Conversely, if $\widetilde{B}:\ Z_0\to Z_0$, then $(\widetilde{B}z)_m=0$ for $m\notin \GD_0$, whenever $z\in Z_0$. Taking $(z_i)=(\delta_{i,j})$, where $j\in \GD_0$ is a fixed element, we find that $\widetilde{b}_{m,j}=0$ for any $j\in \GD_0$. The condition $m\notin \GD_0$ is equivalent to the property that $m\cdot j^{-1}\notin \GD_0$ for any $j\in \GD_0$. Similarly, we check that $\widetilde{b}_{m,j}=0$, if $m\in \GD_0$, $m\cdot j^{-1}\notin \GD_0$.

Obviously, $Z=Z_0\oplus Z_1$, so that $\widetilde{\mathfrak B}$ is isomorphic to the direct sum of the algebras of all linear operators on $Z_0$ and $Z_1$ respectively.. Each of the latter algebras is a type I Baer ring, which implies \cite{Kap} the same property of $\widetilde{\mathfrak B}.\qquad \blacksquare$

\end{document}